\documentclass[12pt,reqno]{amsart}
\usepackage[letterpaper,margin=1.15in]{geometry}
\usepackage{amssymb,mathtools}
\usepackage{microtype}
\usepackage{xcolor}
\usepackage{tikz}
\usepackage{placeins}
\usepackage{accsupp}
\usetikzlibrary{arrows.meta,decorations.pathreplacing}
\usepackage{hyperref}

\definecolor{linkblue}{RGB}{0,85,170}
\definecolor{couplingorange}{RGB}{190,85,0}
\hypersetup{
  colorlinks=true,
  linkcolor=linkblue,
  citecolor=linkblue,
  urlcolor=linkblue,
  pdftitle={Brunn-Minkowski inequality for Riesz capacity in one dimension},
  pdfauthor={Carrie Clark and Richard S. Laugesen},
  pdfsubject={Brunn-Minkowski inequalities for arbitrary compact subsets of the line},
  pdfkeywords={Riesz capacity, logarithmic capacity, Brunn-Minkowski inequality, quantile coupling}
}

\newtheorem{theorem}{Theorem}
\newtheorem{lemma}[theorem]{Lemma}
\theoremstyle{definition}
\newtheorem{example}[theorem]{Example}
\theoremstyle{plain}
\newcommand{\R}{\mathbb R}
\newcommand{\PP}{\mathcal P}
\DeclareMathOperator{\RCap}{Cap}
\DeclareMathOperator{\diam}{diam}
\DeclareMathOperator{\supp}{supp}
\newcommand{\doi}[1]{\href{https://doi.org/#1}{doi:\nolinkurl{#1}}}

\title[Brunn--Minkowski for Riesz capacity]{Brunn--Minkowski inequality for Riesz capacity in one dimension}
\author{Carrie Clark}
\address{Kitchener, Ontario, Canada}
\email{carrieclark435@gmail.com}
\author{Richard S. Laugesen}
\address{University of Illinois, Dept. of Mathematics, Urbana IL 61801, USA}
\email{Laugesen@illinois.edu}
\date{October 24, 2026}

\begin{document}

\begin{abstract}
The Brunn--Minkowski inequality for Riesz $p$-capacity is proved on the line
for every $p<1$ and all compact sets, including disconnected sets.
For $-1<p<1$, equality holds exactly for a singleton summand or else two
intervals, provided at least one capacity is positive.
The proof couples equilibrium measures in increasing order.
\end{abstract}

\maketitle

\vspace{-12pt}

\section{\bf Introduction and main results}\label{sec:introduction}

The classical Brunn--Minkowski inequality in one dimension asserts concavity of 
Lebesgue measure under set addition. We extend this
result to Riesz capacities, which form a one-parameter family of functionals including diameter and logarithmic capacity and also Lebesgue measure as a limiting case. The Riesz parameter $p$ controls the pairwise interaction strength in the kernel $|x-y|^{-p}$, giving greater relative weight to short distances
when $p$ is larger. The logarithmic kernel arises at $p=0$. 

We prove the capacitary Brunn--Minkowski inequality for every $p<1$ and all nonempty compact sets. Equality is characterized precisely. We construct a trial measure on the Minkowski sum by coupling the summands' equilibrium measures so that paired increments always have the same sign. 

Write $\RCap_p$ for Riesz capacity with parameter $p$ (as defined below) and 
\[
K+L=\{x+y:x\in K,\ y\in L\}
\]
for the Minkowski sum of sets $K$ and $L$ in the real line.

\begin{theorem}[Capacitary Brunn--Minkowski on the line]\label{thm:BM}
If $p<1$ then 
\begin{equation}\label{eq:sumBM}
 \RCap_p(K+L)\ge\RCap_p(K)+\RCap_p(L)
\end{equation}
for all nonempty compact sets $K,L\subset\R$.
If $p\le-1$ then equality holds in \eqref{eq:sumBM}.
\end{theorem}

By homogeneity, inequality \eqref{eq:sumBM} for all $K,L$ is equivalent to the concavity formula
\begin{equation}\label{eq:BM}
 \RCap_p\bigl((1-t)K+tL\bigr)
 \ge (1-t)\RCap_p(K)+t\RCap_p(L),\qquad 0\le t\le1,
\end{equation}
for all nonempty compact sets $K,L\subset\R$.

The theorem is proved in Section~\ref{sec:proof}. The new cases are
$-1<p<0$ and $0<p<1$. The logarithmic case $p=0$ follows from Akman's
planar theorem \cite{A26}, as discussed in Section~\ref{sec:literature},
while $p\le-1$ follows from the known diameter formula
\cite[Theorem 1.1(b)]{CLproperties}. Equality is characterized by the following theorem,
proved in Section~\ref{sec:equality}.

\begin{theorem}[Equality conditions]\label{thm:equality-combined}
Fix $-1<p<1$. If $K,L\subset\R$ are nonempty compact sets 
with at least one having positive $p$-capacity when $p\ge0$,
then equality holds in \eqref{eq:sumBM} if and only if
one summand is a singleton or both are closed intervals.
\end{theorem}

The same equality characterization holds in \eqref{eq:BM} for $0<t<1$,
by applying the theorem to $(1-t)K$ and $tL$.

We do not attempt to characterize equality when both $K$ and $L$ have capacity zero. Their sum $K+L$ might or might not have capacity zero. Example~\ref{ex:zero-summands} gives a construction where the sum has positive capacity. 

\subsection*{Definitions}
Given a nonempty compact set $K\subset\R$, denote its collection of Borel
probability measures by $\PP(K)$. Following Clark--Laugesen \cite{CLproperties},
define Riesz capacity by 
\[
 \RCap_p(K)=
 \max_{\mu\in\PP(K)}
 \left(\iint_{K\times K}|x-y|^{-p}\,d\mu(x)\,d\mu(y)\right)^{\!-1/p}, \qquad p \neq 0 ,
\]
under the convention that $\infty^{-1/p}=0$ when $p>0$.  
The maximum can be taken inside the root when $p<0$, but when $p>0$, moving the extremum inside changes it to a minimum. At $p=0$ we use logarithmic capacity:
\[
 \RCap_0(K)=\max_{\mu\in\PP(K)}
 \exp\left(
 - \iint_{K\times K}\log \frac {1}{|x-y|} \,d\mu(x)\,d\mu(y)
 \right),
\]
under the convention that $\exp(-\infty)=0$. 

These maxima are attained by weak-$*$ compactness of $\PP(K)$,
since the functionals being maximized are continuous with respect
to weak-$*$ convergence when $p<0$ and upper semicontinuous
when $p\ge0$. Measures attaining the maxima are called
equilibrium measures, and are interpreted physically as distributions of positive, repelling charges that succeed in ``best spreading out'' the charge over the conductor $K$. 

The double integral in the definition of
$\RCap_p(K)$ is called the $p$-energy of $\mu$, when $p\ne0$, and is called the logarithmic energy when $p=0$. Thus equilibrium measures maximize energy when $p<0$
and minimize it when $p\ge0$.

Capacity is translation invariant, and dilation homogeneous of degree $1$, satisfying 
\[
 \RCap_p(aK+b)=a\RCap_p(K),\qquad a>0,\quad b\in\R.
\] 
For $p<0$ the capacity vanishes exactly on singletons,
since any probability measure that assigned positive charge to two distinct
points would generate a positive double integral.

\subsection*{Exponent range, and endpoint cases for diameter and length}

The interesting parameter range in Theorem~\ref{thm:BM} is $-1<p<1$.
Below the lower endpoint, for $p\le-1$, capacity is simply proportional to diameter with 
$\RCap_p(M)=2^{1/p}\diam(M)$ by
\cite[Theorem 1.1(b)]{CLproperties}. Diameter is obviously linear with respect to 
Minkowski addition on the line.

At the upper endpoint of the range, as $p \uparrow 1$, Theorem~\ref{thm:BM} recovers the classical
Brunn--Minkowski inequality for length, as follows. Writing $m_1$ for Lebesgue
measure on the line we have 
\[
m_1(M) =  \lim_{p\uparrow1} \left(\frac{2}{1-p}\right)^{\!\!1/p}\RCap_p(M)
\]
for every nonempty compact $M\subset\R$, by the one dimensional case of our \cite[Corollary 1.2]{CLproperties}. Thus multiplying \eqref{eq:sumBM} by $(2/(1-p))^{1/p}$ and
letting $p\uparrow1$ yields Brunn--Minkowski for length, 
\[
 m_1(K+L)\ge m_1(K)+m_1(L).
\]

\subsection*{Future directions}

Other authors may wish to seek quantitative versions of
Theorem~\ref{thm:equality-combined}, relating a
capacitary deficit in Brunn--Minkowski to some kind of notion of set distance away from the equality
configurations. 

\subsection*{Structure of the paper}

Section~\ref{sec:literature} situates the results in the literature.
Section~\ref{sec:proof} proves Theorem~\ref{thm:BM} and Section~\ref{sec:equality} proves the equality characterizations in Theorem~\ref{thm:equality-combined}.

\section{\bf Background and prior work}\label{sec:literature}

The classic source for Riesz capacity with $p>0$ is the monograph of Landkof \cite{L72}. His capacity notation and normalization differ from ours. The more recent text by Borodachov, Hardin and Saff \cite[Chapter 4]{BHS19} develops the theory for both positive and negative $p$, again with a differing normalization. Our normalization with the exponent $-1/p$ follows
Clark--Laugesen \cite{CLproperties}. Appendix A in that paper collects known formulas for the capacity of an interval. Also, the diameter formula \cite[Theorem 1.1(b)]{CLproperties} immediately gives
Theorem~\ref{thm:BM} when $p\le-1$.

For logarithmic capacity in the plane, Pommerenke \cite[Satz 2]{P59} proved
Brunn--Minkowski for connected compact sets; see also
\cite[Theorem 11.2(b), p.{\,}336]{P75} for a different proof of the capacity
comparison used in his argument. Borell
\cite[Example 7.3]{B84} later gave a Green function proof for convex sets. Regrettably, when specialized to subsets of the line these results tell us
only about intervals. Akman's recent preprint
\cite[Theorem B(e)]{A26} proves the logarithmic Brunn--Minkowski 
inequality for all compact planar sets. Restricting this inequality to
the line gives the logarithmic case of Theorem~\ref{thm:BM}.
The logarithmic equality cases in Theorem~\ref{thm:equality-combined}
follow from \cite[Theorem A(d) and Remark 6.19(a)]{A26}: a compact
subset of a line has no bounded complementary component in the plane,
and so its filled set is the set itself; the remark handles the case when
exactly one capacity is zero.

Going beyond the one dimensional scope of the current paper to $\R^n, n\ge3$,
Borell \cite{B83} proved Brunn--Minkowski for Newtonian capacity
($p=n-2$) on convex sets. Akman \cite[Theorem A(a)]{A26} extends the
Newtonian result greatly to arbitrary compact sets.
For $n\ge2$, applying the Newtonian theorem in $\R^{n+1}$ to compact
sets in the hyperplane $\R^n$ also gives the Riesz case $p=n-1$,
since the kernel is $|x-y|^{-(n-1)}$. For convex bodies, that case was
treated by Novaga--Ruffini \cite[Theorem 1.1(ii)]{NR15}.
Thus in dimensions $n\ge2$, Brunn--Minkowski is known on 
arbitrary compact sets for Riesz exponents $p=n-2$ and $p=n-1$.

Akman's paper treats not the family of Riesz capacities but rather a family of variational capacities defined by nonlinear Dirichlet integrals and Green functions. These variational capacities overlap with the family of Riesz capacities only at the planar logarithmic and higher dimensional Newtonian cases. 

Some other results for Riesz capacity in one dimension are summarized in Clark--Laugesen \cite[Figure 1]{CLratios}, for the problem of maximizing the ratio $\RCap_q(K)/\RCap_p(K)$ on compact subsets $K \subset \R$. The proved or conjectured maximizers are intervals or two-point sets, 
in different parameter regions. A technical report of Liu \cite[Corollary 1.3]{L26} proves two-point
maximality throughout $q<p<0$.

In a companion paper \cite{CLhigher}, we prove Brunn--Minkowski inequalities for Riesz capacity in higher dimensions. The methods are quite different, because the ordering of the line that underlies the current paper is unavailable in higher dimensions.

\section{\bf Proof of the inequality, Theorem~\ref{thm:BM}}\label{sec:proof}

The key is to pair points in increasing order, so that the distance between two sums is exactly the weighted sum of the corresponding distances, in formula~\eqref{eq:increments} below. Convexity and concavity considerations then provide a trial measure with the required capacity bound. 

When $p\le-1$, equality in \eqref{eq:sumBM} follows immediately from the additivity of diameter, $\diam(K+L)=\diam(K)+\diam(L)$, and the diameter formula $\RCap_p(M)=2^{1/p}\diam(M)$ that we recalled in Section~\ref{sec:introduction}. 

\subsection*{Normalization}
Assume $-1<p<1$. We prove \eqref{eq:sumBM}. 
If either $K$ or $L$ has capacity zero then \eqref{eq:sumBM} follows because $K+L$
contains a translate of each summand. Thus we may assume 
$c=\RCap_p(K)>0$ and $d=\RCap_p(L)>0$. Let 
\[
 \widetilde{K}=c^{-1}K,\qquad \widetilde{L}=d^{-1}L,\qquad t=\frac{d}{c+d} \in (0,1).
\]
These normalized sets have capacity $1$, and
\[
 K+L=(c+d)\bigl((1-t)\widetilde{K}+t\widetilde{L}\bigr).
\]
By homogeneity, the proof of \eqref{eq:sumBM} reduces to showing that the convex combination on the right side has capacity at least $1$. 

Thus we may suppose $\RCap_p(K)=\RCap_p(L)=1$ and $0<t<1$. The goal is to show 
\[
\RCap_p \bigl((1-t)K+t L \bigr) \geq 1 .
\]

\subsection*{Quantile coupling}

Choose equilibrium measures $\mu,\nu$ on $K,L$ respectively. To pair their points in increasing
order, consider the increasing quantile maps $f,g$ on the interval $(0,1)$, defined by 
\[
 f(u)=\min\{x\in\R:\mu((-\infty,x])\ge u\},\qquad 0<u<1,
\]
and similarly for $g$ with $\nu$. Note $f,g$ are bounded and left continuous. These maps push forward Lebesgue measure on
$(0,1)$ to $\mu,\nu$, respectively. Put $h=(1-t)f+t g$
and let $\eta$ be the pushforward under $h$ of Lebesgue measure on $(0,1)$, so that $\eta$ is a Borel probability measure on $\R$. 
Thus $\eta$ is the displacement interpolation between 
$\mu$ and $\nu$, cf.\ McCann \cite[Definition 1.1 and Remark 1.8]{McC97}.

We show $\eta$ is supported on the convex combination $C=(1-t)K+t L$. For this, observe that the pushforward property and $\mu(K)=\nu(L)=1$ give
$f(u)\in K$ and $g(u)\in L$ for almost every $u \in (0,1)$.
Hence $h(u)\in C$ a.e., and so $\eta(C)=1$.
The set $C$ is compact, being the continuous image of $K\times L$
under $(x,y)\mapsto(1-t)x+t y$, and therefore $\supp\eta\subset C$.

\begin{figure}[!t]
\centering
\BeginAccSupp{method=pdfstringdef,unicode,Alt={Schematic of ordered quantile coupling. Three number lines show two ordered points from K, from L and from their Minkowski average. For u less than v, the left points are f(u), g(u) and h(u), while the right points are f(v), g(v) and h(v). Filled markers and solid guides identify the u samples. Open markers and dashed guides identify the v samples. The distances D0 and D1 combine to D equal to (1 minus t) times D0 plus t times D1, because order is preserved.},ActualText={Schematic of ordered quantile coupling. Three number lines show two ordered points from K, from L and from their Minkowski average. For u less than v, the left points are f(u), g(u) and h(u), while the right points are f(v), g(v) and h(v). Filled markers and solid guides identify the u samples. Open markers and dashed guides identify the v samples. The distances D0 and D1 combine to D equal to (1 minus t) times D0 plus t times D1, because order is preserved.}}
\begin{tikzpicture}[x=0.92cm,y=1cm,font=\small,
  axis/.style={-{Stealth[length=1.8mm]},line width=0.4pt,black!55},
  guide/.style={line width=0.65pt},
  pointlabel/.style={below=3pt,fill=white,inner sep=1.5pt},
  distancebrace/.style={decorate,decoration={brace,amplitude=4pt},line width=0.5pt}]
  \draw[guide,linkblue] (1,4.6)--(2.5,2.3)--(1.75,0);
  \draw[guide,couplingorange,densely dashed] (6,4.6)--(9,2.3)--(7.5,0);
  \draw[axis] (-0.1,4.6)--(10.2,4.6);
  \draw[axis] (-0.1,2.3)--(10.2,2.3);
  \draw[axis] (-0.1,0)--(10.2,0);
  \node[anchor=east] at (-0.45,4.6) {$K$};
  \node[anchor=east] at (-0.45,2.3) {$L$};
  \node[anchor=east] at (-0.45,0) {$(1-t)K+t L$};
  \fill[linkblue] (1,4.6) circle (2pt);
  \draw[couplingorange,fill=white,line width=0.9pt] (6,4.6) circle (2.3pt);
  \fill[linkblue] (2.5,2.3) circle (2pt);
  \draw[couplingorange,fill=white,line width=0.9pt] (9,2.3) circle (2.3pt);
  \fill[linkblue] (1.75,0) circle (2pt);
  \draw[couplingorange,fill=white,line width=0.9pt] (7.5,0) circle (2.3pt);
  \node[pointlabel,text=linkblue] at (1,4.6) {$f(u)$};
  \node[pointlabel,text=couplingorange] at (6,4.6) {$f(v)$};
  \node[pointlabel,text=linkblue] at (2.5,2.3) {$g(u)$};
  \node[pointlabel,text=couplingorange] at (9,2.3) {$g(v)$};
  \node[pointlabel,text=linkblue] at (1.75,0) {$h(u)$};
  \node[pointlabel,text=couplingorange] at (7.5,0) {$h(v)$};
  \draw[distancebrace] (1,5.05)--(6,5.05)
    node[midway,above=6pt,fill=white,inner sep=1pt] {$D_0$};
  \draw[distancebrace] (2.5,2.75)--(9,2.75)
    node[midway,above=6pt,fill=white,inner sep=1pt] {$D_1$};
  \draw[distancebrace] (1.75,0.45)--(7.5,0.45)
    node[midway,above=6pt,fill=white,inner sep=1pt]
       {$D=(1-t)D_0+t D_1$};
\end{tikzpicture}
\EndAccSupp{}
\caption{The quantile maps preserve order and so the paired
increments add without cancellation, in the weighted sum.}
\label{fig:ordered-coupling}
\end{figure}
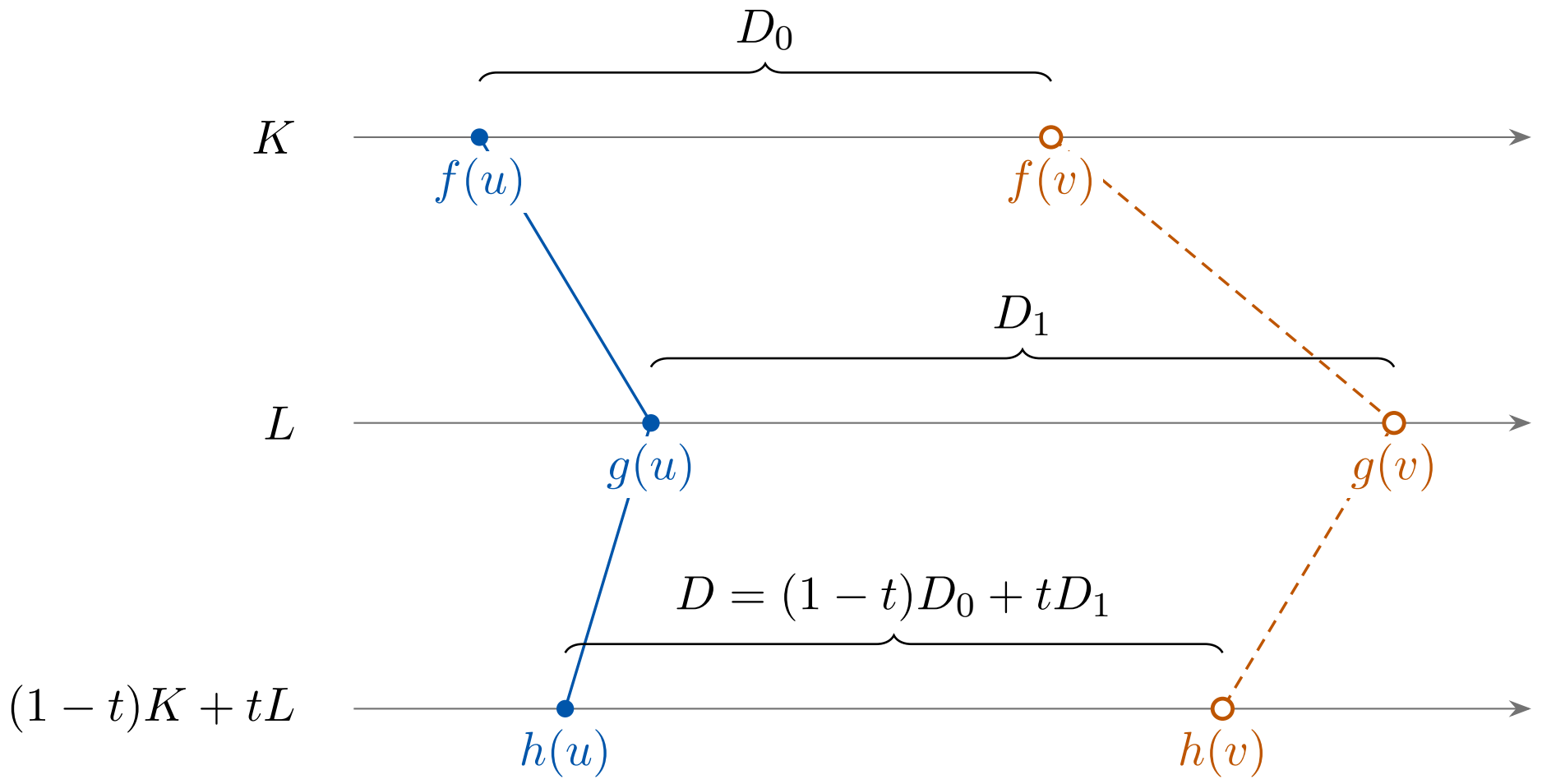

Put $D_0(u,v)=|f(u)-f(v)|$ and $D_1(u,v)=|g(u)-g(v)|$.
Since $f,g$ are increasing, their paired increments have the same sign,
and so
\begin{equation}\label{eq:increments}
 |h(u)-h(v)|
 =(1-t)D_0(u,v)+t D_1(u,v)=:D(u,v).
\end{equation}
Figure~\ref{fig:ordered-coupling} illustrates this distance relation.

For $p\ne0$, the capacity equaling $1$ for $K,L$ means their energies equal $1$, so that 
\[
 \iint_{(0,1)^2}D_j^{-p}\,du\,dv=1,\qquad j=0,1.
\]
By convexity of $r\mapsto r^{-p}$ when $0<p<1$ (the displacement convexity argument of
\cite[Proposition 2.7, Remark 3.4]{CFP12}) and concavity of $r\mapsto r^{-p}$ when $-1<p<0$, we deduce
\begin{equation} \label{eq:concaveconvex}
 \iint_{(0,1)^2}D^{-p}\,du\,dv
 \begin{cases}
 \ge1,&-1<p<0,\\
 \le1,&0<p<1.
 \end{cases}
\end{equation}
In each case the trial measure $\eta$ on $C$ yields capacity greater than or equal to $1$.

For $p=0$, the energy integral of $\log(1/D_j)$ equals $0$ for $j=0,1$.
Using convexity of $\log(1/r)$ (which is the logarithmic displacement
convexity argument of \cite{Bl04,CFP12}) implies that the trial measure $\eta$
has energy $\iint_{(0,1)^2}\log(1/D) \,du\,dv \le 0$, and hence $C$ has
logarithmic capacity at least $1$.

\emph{Note.} The convexity and concavity are strict. Equality therefore requires
$D_0=D_1$ almost everywhere with respect to Lebesgue measure
$du\,dv$ on $(0,1)^2$ and requires $\eta$ to be an equilibrium measure. These conditions will determine
the equality cases in Section~\ref{sec:equality}.

\FloatBarrier
\section{\bf Proof of equality characterizations, Theorem~\ref{thm:equality-combined}}\label{sec:equality}

Equality presents two rigidity questions. First, the convexity and concavity 
inequalities in Section~\ref{sec:proof} force the equilibrium measures
to be positively homothetic, as we show below. We must then recover the
sets from their measures and use equality once more to rule out gaps.
The following lemma supplies the needed support information.

Fan \cite[Theorem 4.3]{F26} proves positivity of every equilibrium
charge for finite subsets of the line when $-1<p<0$.
The negative parameter part of the next lemma extends her argument
to arbitrary compact subsets of the line and gives strict inclusion
monotonicity.

Given a probability measure
$\mu$ on a nonempty compact set $K\subset\R$, write the support, potential and energy as 
\[
 S=\supp\mu,\qquad
 U(x)=\int_K|x-y|^{-p}\,d\mu(y)\quad(x\in\R),\qquad
 E=\int_K U(x)\,d\mu(x).
\]
When $p=0$, replace the kernel with $\log(1/|x-y|)$.

\begin{lemma}[Equilibrium support]\label{lem:equilibrium-support}
Suppose $-1<p<1$. If $K\subset\R$ is nonempty and compact and $\mu$ is an
equilibrium measure on $K$ with finite energy, then every probability measure on
$K$ with finite energy is concentrated in the support $S$ of $\mu$. Moreover,
\begin{equation} \label{eq:pointwisebounds}
 \begin{split}
 U&\ge E\text{ on }\R \text{\ and \ } U>E\text{ on }\R\setminus S \quad 
 \text{if }-1<p<0,\\
 U&\le E\text{ on }\R \text{\ and \ } U<E\text{ on }\R\setminus S \quad 
 \text{if } 0 \leq p < 1.
 \end{split}
\end{equation}
Further, if $-1<p<0$ then $S=K$ and
\[
 \RCap_p(K)<\RCap_p(M)
\]
whenever $M \subset \R$ is compact and $K\subsetneqq M$.
\end{lemma}

The finite energy hypothesis on $\mu$ is automatic when $p<0$, and when $p\ge0$ it is equivalent
to $\RCap_p(K)>0$.

\begin{proof}
When $p<0$, the book of Borodachov, Hardin and Saff
\cite[Definitions 2.2.1 and 4.1.4]{BHS19} uses $-|x-y|^{-p}$, which is the negative of our kernel, 
and hence they define equilibrium measures by energy minimization rather than maximization. Their potential and energy are $-U$ and $-E$ in our notation, when $p<0$. The equilibrium inequalities in
\cite[Lemma 4.2.13 and Theorems 4.2.14--4.2.15]{BHS19}
are stated for infinite compact sets, but their proofs apply
unchanged to finite compact sets. We use these results below
in this generality. 

Now we combine results from \cite{BHS19} to obtain the lemma. By \cite[Theorem 4.2.15]{BHS19},
\[
 U\ge E\quad\text{on }S\quad(p<0),\qquad
 U\le E\quad\text{on }S\quad(p\ge0).
\]
The First Maximum Principle \cite[Theorem 4.5.11]{BHS19} extends these
bounds to $\R$, applied in dimension $1$ to $-U$ when $p<0$ and to $U$
when $p\ge0$.
On each component of $\R\setminus S$, the potential $U$ is strictly concave when $p<0$
and strictly convex when $p\ge0$, and so it cannot attain the corresponding
lower or upper bound at an interior point. This argument proves \eqref{eq:pointwisebounds}.

Now suppose $\sigma$ is a probability measure on $K$ with finite energy. Applying \cite[Theorem 4.2.14]{BHS19} to the equilibrium measure \(\mu\) and using the inequalities from \eqref{eq:pointwisebounds} gives that $U=E$ on $K$ outside an exceptional set (when $p \geq 0$) that is a countable union of compact sets having zero Riesz $p$-capacity. That exceptional set has $\sigma$-measure zero by \cite[Lemma 4.2.13]{BHS19}, and so in all cases $U= E$ $\sigma$-almost everywhere. 
Combining this fact with the strict bounds in \eqref{eq:pointwisebounds} outside $S$ gives that 
$\sigma(K\setminus S)=0$.
When $p<0$, applying this conclusion with $\sigma$ a Dirac measure shows that $K \setminus S$ is empty, that is, $S=K$. 

Finally, if $p<0$, $K\subset M\subset\R$ and the capacities of $K$ and $M$ agree, then $\mu$ is an equilibrium measure for both $K$ and $M$. Hence $K=S=M$, proving the contrapositive of strict inclusion monotonicity.
\end{proof}

\subsection*{\texorpdfstring
  {Proof of Theorem~\ref{thm:equality-combined} ``$\Longleftarrow$''}
  {Proof of Theorem~\ref{thm:equality-combined}: sufficiency}}\ 
  
Translation invariance gives equality in \eqref{eq:sumBM} when one set is a singleton.
For two intervals, the lengths add and capacity is proportional to
length by homogeneity, again giving equality.

\subsection*{\texorpdfstring
  {Proof of Theorem~\ref{thm:equality-combined}
   ``$\Longrightarrow$'' for $-1<p<0$}
  {Proof of Theorem~\ref{thm:equality-combined}:
   necessity for -1 < p < 0}}\ 
   
Suppose neither set is a singleton, so that both sets have positive capacity, and that equality holds in
\eqref{eq:sumBM}. Using the normalization and notation from the proof of
Theorem~\ref{thm:BM}, we relabel the normalized sets as $K,L$.
Then $0<t<1$ and
$\RCap_p(K)=\RCap_p(L)=\RCap_p(C)=1$.
Examining that proof, we see equality holds in \eqref{eq:concaveconvex}. Hence strict concavity of $r\mapsto r^{-p}$ implies
\[
|f(u)-f(v)| =|g(u)-g(v)|
 \quad\text{for almost every }(u,v)\in(0,1)^2.
\]
The increments have the same sign, since $f,g$ are increasing, and so 
\begin{equation}\label{eq:constant-offset}
 g(u)-f(u)=g(v)-f(v)
 \quad\text{for almost every }(u,v)\in(0,1)^2.
\end{equation}
Integrating \eqref{eq:constant-offset} with
respect to $v$ gives $g(u)-f(u)=\text{constant}$ for almost every $u\in(0,1)$. Left continuity extends this identity to every $u\in(0,1)$. We may suppose the constant is zero, by translating $L$. 
Then $f=g=h$ pointwise and $\mu=\nu=\eta$. Lemma~\ref{lem:equilibrium-support} gives $K=L=\supp\mu$.
Both the energy of $\mu$ and $\RCap_p(C)^{-p}$ equal $1$, and so
$\mu$ is also an equilibrium measure on $C$.
Applying Lemma~\ref{lem:equilibrium-support} again gives $C=K$.
Since $K=L$, the definition of $C$ now yields
\begin{equation}\label{eq:weighted-sum-closure}
K = C = (1-t)K+tK.
\end{equation}

If $K$ were not an interval, there would be a gap $(\alpha,\beta)$ of $K$
between its extreme points, with $\alpha,\beta\in K$. The point $(1-t)\alpha+t\beta$ lies in that gap, whereas \eqref{eq:weighted-sum-closure} requires it to
belong to $K$. This contradiction shows that $K$ is an interval, and so
$L$ is also an interval. The original $K$ and $L$ are therefore intervals too. 

\subsection*{\texorpdfstring
  {Proof of Theorem~\ref{thm:equality-combined}
   ``$\Longrightarrow$'' for $0\le p<1$}
  {Proof of Theorem~\ref{thm:equality-combined}:
   necessity for 0 <= p < 1}}\ 
   
Equality of the equilibrium measures alone need not
identify the sets, since a measure can ignore parts of capacity zero.
We handle those parts using the support conclusion of
Lemma~\ref{lem:equilibrium-support}.

First suppose exactly one of the capacities is zero, say
 $\RCap_p(K)>0$ and $\RCap_p(L)=0$, and that equality holds in \eqref{eq:sumBM}. Translate $L$ so that $0\in L$. An equilibrium measure $\mu$ on
$K$ is then an equilibrium measure on $K+L$, since their capacities
are equal. Write $S=\supp\mu$. For each $y\in L$, the translate of
$\mu$ by $y$ has finite energy and is supported on $K+L$. Applying 
Lemma~\ref{lem:equilibrium-support} to $K+L$ therefore gives $S+y\subset S$. By considering the maximum and minimum points of $S$, one infers that $y=0$. Thus $L=\{0\}$ is a singleton.

Now suppose both capacities are positive and equality holds in
\eqref{eq:sumBM}. Argue as in the negative parameter case, using strict convexity
of $r\mapsto r^{-p}$ for $p>0$ and of $\log (1/r)$ for $p=0$.
After normalization and translation, $f=g=h$ and $\mu=\nu=\eta$ is
a common equilibrium measure on $K,L$ and $C=(1-t)K+tL$, where $0<t<1$.
Lemma~\ref{lem:equilibrium-support} implies that every finite
energy probability measure on $C$ is concentrated on
$S=\supp\mu\subset K\cap L$. Note $S$ is not a singleton, since $\mu$ has finite energy. 

If $S$ has a gap $(\alpha,\beta)$ between its extreme points, choose a
closed neighborhood $Q$ of $\alpha$ small enough that
\[
 (1-t)Q+t\beta\subset(\alpha,\beta).
\]
Since $\alpha\in S$, one has $\mu(Q)>0$. The normalized restriction of
$\mu$ to $Q$ has finite energy. The pushforward of the measure under $x\mapsto(1-t)x+t\beta$ also has
finite energy and is supported in $C\cap(\alpha,\beta)$, contradicting
Lemma~\ref{lem:equilibrium-support}. Thus $S$ is an interval, of positive
length because $\mu$ has finite energy.

If $x\in K\setminus S$, then $x$ lies outside this interval.
The interval $(1-t)x+tS\subset C$ has a portion of positive length
outside $S$. Normalized Lebesgue measure on a closed subinterval of
that portion has finite energy and thus contradicts
Lemma~\ref{lem:equilibrium-support}. Hence $K=S$, and similarly $L=S$. Undoing the normalization and translation shows that the
original summands are intervals, giving the equality condition for \eqref{eq:sumBM}. 

\begin{example}[Zero-capacity sets can have positive-capacity sum]\label{ex:zero-summands}
Let $A\subset[0,1]$ be the middle thirds Cantor set and suppose
$\log 2/\log 3<p<1$. We claim
\[
 \RCap_p(A)=0,\qquad A+A=[0,2],\qquad \RCap_p(A+A)>0, 
\]
so that the Brunn--Minkowski inequality is strict in this case.

Indeed, the Hausdorff dimension of $A$ is $\log 2/\log 3<p$ and so
$\RCap_p(A)=0$ by \cite[Theorem 4.3.1]{BHS19}.
The classical identity $A+A=[0,2]$ is recorded in
\cite[Exercise 3.3.7]{Ab15}. That interval has positive capacity, for example by taking normalized Lebesgue measure as a trial measure. 

For an example that works for every $0\le p<1$, note that Gorodetski--Luna--Wood
\cite[Theorem 1]{GLW26} construct a different Cantor set $B\subset[0,1]$
with $\RCap_0(B)=0$ and $B+B=[0,2]$. Capacity is decreasing with respect
to $p$ by \cite[Theorem 1.1(a)]{CLproperties}, and so $\RCap_p(B)=0$
for every $p>0$. Thus $B+B$ gives strict inequality in
\eqref{eq:sumBM} for each $p \in[0,1)$. 
\end{example}

\section*{\bf Acknowledgments}

Laugesen's research was supported by grants from the National Science Foundation (\#2246537), Simons Foundation (\#964018) and the University of Illinois Campus Research Board (RB26003). This work was supported in part by a grant of access to OpenAI models
through the \emph{ChatGPT for Academic Researchers} program.

\section*{\bf AI usage statement}

ChatGPT 6.1-Sol helped to develop the mathematical arguments, check
calculations, and search the literature, and it assisted with drafting
and revising the manuscript.


\begin{thebibliography}{99}

\bibitem{Ab15}
S. Abbott,
Understanding Analysis, second edition,
Undergraduate Texts in Mathematics, Springer, New York, 2015.
\doi{10.1007/978-1-4939-2712-8}

\bibitem{A26}
M. Akman,
\emph{Brunn--Minkowski inequality for $p$-capacity of compact sets},
arXiv:2610.28421 (2026),
\url{https://arxiv.org/abs/2610.28421}.

\bibitem{Bl04}
G. Blower,
\emph{Displacement convexity for the generalized orthogonal ensemble},
J. Stat. Phys. 116 (2004), 1359--1387.
\doi{10.1023/B:JOSS.0000041742.86859.cd}

\bibitem{B83}
C. Borell,
\emph{Capacitary inequalities of the Brunn--Minkowski type},
Math. Ann. 263 (1983), 179--184.
\doi{10.1007/BF01456879}

\bibitem{B84}
C. Borell,
\emph{Hitting probabilities of killed Brownian motion: a study on geometric
regularity},
Ann. Sci. \'{E}cole Norm. Sup. (4) 17 (1984), 451--467.
\doi{10.24033/asens.1480}

\bibitem{BHS19}
S. V. Borodachov, D. P. Hardin and E. B. Saff,
Discrete Energy on Rectifiable Sets,
Springer Monographs in Mathematics, Springer, New York, 2019.
\doi{10.1007/978-0-387-84808-2}

\bibitem{CFP12}
J. A. Carrillo, L. C. F. Ferreira and J. C. Precioso,
\emph{A mass-transportation approach to a one dimensional fluid mechanics
model with nonlocal velocity},
Adv. Math. 231 (2012), 306--327.
\doi{10.1016/j.aim.2012.03.036}

\bibitem{CLratios}
C. Clark and R. S. Laugesen,
\emph{Maximizing Riesz capacity ratios: conjectures and theorems},
SIAM J. Math. Anal. 57 (2025), 5683--5717.
\doi{10.1137/24M171992X}

\bibitem{CLproperties}
C. Clark and R. S. Laugesen,
\emph{Riesz capacity: monotonicity, continuity, diameter and volume},
Anal. Math. Phys. 15 (2025), Paper No. 4.
\doi{10.1007/s13324-024-01000-2}

\bibitem{CLhigher}
C. Clark and R. S. Laugesen,
\emph{Brunn--Minkowski inequalities for Riesz capacity in all dimensions},
arXiv:2611.19258 (2026),
\url{https://arxiv.org/abs/2611.19258}.

\bibitem{F26}
Q. Fan,
\emph{Riesz capacity ratios with negative exponents},
arXiv:2609.11186 (2026),
\doi{10.48550/arXiv.2609.11186}.

\bibitem{GLW26}
A. Gorodetski, A. Luna and W. Wood,
\emph{On sums of ``small'' Cantor sets},
arXiv:2610.17005 (2026),
\doi{10.48550/arXiv.2610.17005}.

\bibitem{L72}
N. S. Landkof,
Foundations of Modern Potential Theory,
translated from the Russian by A. P. Doohovskoy.
Die Grundlehren der mathematischen Wissenschaften, Band 180,
Springer--Verlag, New York--Heidelberg, 1972.
\doi{10.1007/978-3-642-65183-0}

\bibitem{L26}
Z. Liu (producer),
\emph{Sharp ratios and stability for negative Riesz capacities},
working paper, Zenodo, version 1 (2026).
\doi{10.5281/zenodo.23211695}

\bibitem{McC97}
R. J. McCann,
\emph{A convexity principle for interacting gases},
Adv. Math. 128 (1997), 153--179.
\doi{10.1006/aima.1997.1634}

\bibitem{NR15}
M. Novaga and B. Ruffini,
\emph{Brunn--Minkowski inequality for the $1$-Riesz capacity and level set
convexity for the $1/2$-Laplacian},
J. Convex Anal. 22 (2015), 1125--1134.
\doi{10.68381/jca22059}

\bibitem{P59}
C. Pommerenke,
\emph{\"{U}ber die Kapazit\"{a}t der Summe von Kontinuen},
Math. Ann. 139 (1959), 127--132.
\doi{10.1007/BF01354870}

\bibitem{P75}
C. Pommerenke,
Univalent Functions,
Studia Mathematica/Mathematische Lehrb\"ucher, Band XXV,
Vandenhoeck \& Ruprecht, G\"ottingen, 1975.

\end{thebibliography}
\end{document}